\documentclass[10pt,a4paper]{amsart}

\usepackage{amssymb,amsmath,amsfonts}

\usepackage{graphicx}

\usepackage{mathpazo}
\usepackage[hyperindex,pageanchor]{hyperref}

\numberwithin{equation}{section}
\newtheorem{theorem}{Theorem}[section]
\newtheorem{lemma}[theorem]{Lemma}
\newtheorem{proposition}[theorem]{Proposition}
\newtheorem{corollary}[theorem]{Corollary}

\theoremstyle{definition}
\newtheorem{definition}[theorem]{Definition}
\theoremstyle{remark}
\newtheorem{remark}[theorem]{Remark}

\newtheorem{acknowledgement}{Acknowledgement}

\newcommand{\Ass}{\operatorname{Ass}}

\newcommand{\grade}{\operatorname{grade}}

\newcommand{\Assh}{\operatorname{Assh}}
\newcommand{\Spec}{\operatorname{Spec}}

\newcommand{\Ht}{\operatorname{ht}}

\newcommand{\Ext}{\operatorname{Ext}}
\newcommand{\Supp}{\operatorname{Supp}}
\newcommand{\CoSupp}{\operatorname{CoSupp}}

\newcommand{\Hom}{\operatorname{Hom}}
\newcommand{\Att}{\operatorname{Att}}
\newcommand{\Ann}{\operatorname{Ann}}
\newcommand{\Rad}{\operatorname{Rad}}

\newcommand{\depth}{\operatorname{depth}}

\newcommand{\Coass}{\operatorname{Coass}}
\newcommand{\coker}{\operatorname{coker}}
\newcommand{\Cosupp}{\operatorname{Cosupp}}

\newcommand{\vpl}{\operatornamewithlimits{\varprojlim}}

\newcommand{\vil}{\operatornamewithlimits{\varinjlim}}

\newcommand{\fm}{\frak{m}}
\newcommand{\fp}{\frak{p}}
\newcommand{\fq}{\frak{q}}
\newcommand{\fa}{\frak{a}}
\newcommand{\fb}{\frak{b}}

\newcommand{\fn}{\frak{n}}

\begin{document}

\author[Eghbali ]{Majid Eghbali  }

\title[Artinian formal local cohomology]
{ On Artinianness of formal local cohomology, colocalization and coassociated primes }

\address{Martin-Luther-Universit\"at Halle-Wittenberg, Institut f\"ur Informatik, D -- 06099 Halle (Saale),
 Germany.-And- Islamic
Azad University Branch Sofian, Sofian, Iran}
\email{m.eghbali@yahoo.com}



\begin{abstract}  This paper at first concerns some criteria on  Artinianness and vanishing of formal local cohomology modules. Then we consider the cosupport and the set of coassociated primes of these modules more precisely.
\end{abstract}

\maketitle

\section{Introduction}

Throughout, $\fa$ is an ideal of a commutative Noetherian ring $R$ and $M$  an $R$-module. Let $V(\fa)$ be the set of prime ideals in $R$ containing $\fa$. For an integer $i$, let $H^{i}_{\fa}(M) \ $ denote the $i$-th local cohomology module of $M$. We have the isomorphism of $H^{i}_{\fa}(M) \ $  to   ${\vil}_n \Ext^{i}_{R}(R/\fa^n, M)$ for every $i \in \mathbb {Z}$, see \cite{Br-Sh} for more details.
%Herzog \cite{5}  introduced
%$$\begin{array}{ll} \ H^{i}_{\fa}(M , N) := {\vil}_n \Ext^{i}_{R}(M/\fa^n M,N), \ i\in  \mathbb{Z}
%\end{array}$$
%for two $R$-modules $M$, $N$ and an ideal $\fa$ as the generalized local cohomology module with respect to $\fa$. For $R=M$ it is equal to the ordinary local cohomology of $N$ with respect ideal $\fa$.

Consider the family of local cohomology modules $\{ H^{i}_{\fm}(M/{\fa}^n M) \}_{n \in \mathbb{N}} \ $. For every $n$ there is a natural homomorphism $H^{i}_{\fm}(M/{\fa}^{n+1} M) \rightarrow H^{i}_{\fm}(M/{\fa}^n M)$ such that the family forms a projective system. The projective limit
$\mathfrak{F}^i_{\fa}(M):={\vpl}_nH^i_{\fm}(M/\fa^n M)$ is called the $i$-th formal local cohomology of $M$ with respect to $\fa$. Formal local cohomology modules were used by Peskine and Szpiro in \cite{P-S} when $R$ is a regular ring in order to solve a conjecture of Hartshorne in prime characteristic. It is noteworthy to mention that if
$U = \Spec(R) \setminus \{\fm\}$ and $( \widehat{U},\mathcal{O}_{\widehat{u}})$ denote the formal completion of $U$ along $V(\fa) \setminus \{\fm \}$
and also $\widehat{\mathcal{F}}$ denotes the $\mathcal{O}_{\widehat{u}}$-sheaf associated to ${\vpl}_n M/\fa^n M$, they have described the formal cohomology modules $H^i(\widehat{U},\mathcal{O}_{\widehat{u}})$ via the isomorphisms $H^i(\widehat{U},\mathcal{O}_{\widehat{u}})\cong \mathfrak{F}^i_{\fa}(M)$, $i \geq 1$. See also   \cite[proposition (2.2)]{O} when $R$ is a Gorenstein ring.

Let $\textbf{x}= \{ x_1,..., x_r \}$ denote a system of elements such that $\fm= \Rad (\textbf{x})$. In  \cite{Sch}, Schenzel has studied formal local cohomology module via following isomorphism
$$\begin{array}{ll} \ {\vpl}_nH^i_{\fm}(M/\fa^n M)\cong H^{i} ({\vpl}_n (\check{C_{\textbf{x}}} \otimes M/\fa^n M))
\end{array}$$
where  $\ \check{C_{\textbf{x}}}\ $ denotes the $\check{C}$ech complex of $R$ with respect to $\textbf{x}$ .

When the local ring $(R,\fm)$ is a quotient of a local Gorenstein ring $(S,\fn)$, we have
$$\begin{array}{ll} \ \mathfrak{F}^i_{\fa}(M)\cong \Hom_{R}(H^{\dim S-i}_{\fa^{'}}(M , S),E),  \ i\in  \mathbb{Z} \ \ \ \ \ \ \ \ \ (1.1)
\end{array}$$
where $E$ denotes the injective hull of $R/\fm$ and $\fa^{'}$ is the preimage of $\fa$ in $S$ (cf. \cite[Remark 3.6]{Sch}).

Important problems concerning local cohomology modules  are vanishing, finiteness and  Artinianness results (see, e.g., \cite{Hu}). In Section 2 we examine the vanishing and Artinianness of  formal local cohomology modules. In the next theorem we give  some criteria for vanishing and
Artinianness of formal local cohomology modules:

\begin{theorem} \label{1.1}  Let  $(R,\fm)$ be a  local ring and $M$ be a finitely generated $R$-module. For given integers $i$ and $t > 0$,  the following statements are equivalent:
\begin{enumerate}
\item [(1)] $\Supp_{\widehat{R}} (\mathfrak{F}^i_{\fa}(M)) \subseteq V(\fm \widehat{R})$ for all $i < t$;
\item[(2)] $\mathfrak{F}^i_{\fa}(M)$ is Artinian for all $i < t$;
\item[(3)] $\Supp_{\widehat{R}} (\mathfrak{F}^i_{\fa}(M)) \subseteq V(\fa \widehat{R})$  for all $i < t$;
\item[(4)] $\fa \subseteq \Rad(\Ann_R (\mathfrak{F}^i_{\fa}(M)))$ for all $i < t$;

Suppose that   $t \leq \depth M$, then the above conditions are equivalent to
 \item[(5)] $\mathfrak{F}^i_{\fa}(M) = 0 $ for all $i < t$;
\end{enumerate}
where $\widehat{R}$ denotes the $\fm$-adic completion of $R$.
\end{theorem}

It should be noted that  it has been shown independently in \cite {Ma} that statements (2) and (4) are equivalent.

Note that as we see in Theorem \ref{1.1}, we have the equivalence between   $\Supp_{\widehat{R}} (\mathfrak{F}^i_{\fa}(M)) \subseteq V(\fa \widehat{R})$  for all $i < t$ and $\fa \subseteq \Rad(\Ann_R (\mathfrak{F}^i_{\fa}(M)))$ for all $i < t$, which is not true in general for an arbitrary module.

In Section 3,  we study the cosupport of formal local cohomology via Richardson's definition of colocalization (cf. Definition \ref{3.1}). We show that when $(R,\fm)$ is a local ring, $M$ is a finite $R$-module and $\mathfrak{F}^i_{\fa}( M)$ is Artinian $(i \in \mathbb {Z})$, then $\CoSupp (\mathfrak{F}^i_{\fa}( M)) \subseteq V(\fa)$ (cf. \ref{3.5}). As a further result we reduce to the case $M=R$ when considering the cosupport of top formal local cohomology modules which is the analogue for formal local cohomology of the result due
to Huneke-Katz-Marley in \cite [Proposition 2.1]{Hu-K-Mar}:

\begin{theorem} \label{1.2} Let $(R, \fm)$ be a  local  ring. Let $M$  be a  finitely generated  $R$-module. Then
\begin{enumerate}
\item[(1)] $\CoSupp (\mathfrak{F}^c_{\fa}( M)) = \CoSupp (\mathfrak{F}^c_{\fa}( R/J))$,
\item[(2)] $\Supp (\mathfrak{F}^c_{\fa}( M)) = \Supp (\mathfrak{F}^c_{\fa}( R/J))$,
\end{enumerate}
where $J$ is $\Ann_{R} M$ and $c:= \dim R/\fa$.
\end{theorem}

For a  representable module $M$,  $\CoSupp M = V(\Ann M)$ (cf. \cite [Theorem 2.7]{R}). It motivates us to see when the cosupport of formal local cohomology module is a closed subset of $\Spec R$ in Zariski topology. For this reason in Section 4 we study the set of coassociated  primes of formal local cohomology more precisely. In this direction  when $(R, \fm)$ is a local ring and $M$ is an $R$-module, the set of  minimal primes in $\CoSupp(M)$ is finite if and only if $\CoSupp(M)$ is a closed subset of $\Spec R$ (lemma \ref{4.2}). Hence, it is enough to ask when the $\Coass M$ is finite. We give affirmative answers to this question in some cases, see Proposition \ref{4.4} and Theorem \ref{1.4} below. It is noteworthy that for a finitely generated module $M$ over a local ring $(R, \fm)$, $\Coass_{\widehat{R}} \mathfrak{F}^{0}_{\fa}(M)$ is finite  since  $\mathfrak{F}^{0}_{\fa}(M)$ is a finite $\widehat{R}$-module (cf. \cite [Lemma 4.1]{Sch}) and $\Coass \mathfrak{F}^{\dim}_{\fa}(M)$ is finite as $ \mathfrak{F}^{\dim M}_{\fa}(M)$ is an Artinian $R$-module (cf. \cite [Lemma 2.2]{A-D} or Proposition \ref{2.1}).

As final results in Section 4, we give the following results for top formal local cohomology modules:

\begin{theorem} \label{1.3} Let $(R,\fm)$ be a local ring of dimension $d > 1$. Let $\mathfrak{F}^{d}_{\fa}(R)=0$. Then:
\begin{enumerate}
\item[(1)] If $\fp \in \Coass \mathfrak{F}^{d-1}_{\fa}(R)$, then it implies that $\dim (R/(\fa,\fp))= d-1$.
\item[(2)] $\Assh(R) \cap \Coass \mathfrak{F}^{d-1}_{\fa}(R)\subseteq \{ \fp \in \Spec R: \dim R/\fp=d, \Rad (\fa+ \fp) \neq \fm \}$.
\item[(3)] If  $ \Coass \mathfrak{F}^{d-1}_{\fa}(R)\subseteq \Assh(R)$, then $\{ \fp \in \Spec R: \dim (R/(\fa,\fp))= d-1 \} \subseteq \Coass \mathfrak{F}^{d-1}_{\fa}(R)$.
\end{enumerate}
\end{theorem}

Next result shows that for a one dimensional ideal $\fa$ of a complete local ring $R$ of dimension $d$,  $\Cosupp \mathfrak{F}^{d-1}_{\fa}(R)$ is closed.

\begin{theorem}  \label{1.4} Let $(R,\fm)$ be a local complete ring of dimension $d$. Let  $\fa$ be an ideal of dimension one. Then
 $$\begin{array}{ll} \  \mathfrak{F}^{d-1}_{\fa}(R)=0, \   when \ \  d >2,
\end{array}$$
in particular $\Coass_R \mathfrak{F}^{d-1}_{\fa}(R)= \emptyset$.
 $$\begin{array}{ll} \ \Coass_R \mathfrak{F}^{d-1}_{\fa}(R) \subseteq \{ \fm \}, \   when\ \  d =1,
\end{array}$$
and in the case $d=2$ we have
 $$\begin{array}{ll} \ \Coass_R \mathfrak{F}^{d-1}_{\fa}(R)= \bigcup^r_{i=1} \Coass_R(\widehat{R_{\fp_i}}) =  \{ \fp_1,..., \fp_r \}  \cup  (\bigcup^s_{j=1} \{\fq_j : R_{\fp_i}/\fq_j R_{\fp_i} \text{ is not complete } \}) ,
\end{array}$$
where $\fp_1 ,... \fp_r$  are minimal prime ideals of  $\fa$ and  $\fq_1 ,... \fq_s$  are minimal prime ideals of  $R$ with $\fq_j \subseteq \fp_i$ for $i \in \{1,...,r \}$.

In particular $\Cosupp \mathfrak{F}^{d-1}_{\fa}(R)$ is closed for all $d >0$.
\end{theorem}

We denote by $\max(R)$ the set of maximal ideals of $R$.

\section{On Artinianness of $\mathfrak{F}^i_{\fa}(M)$}

Important problems concerning local cohomology modules are vanishing,
finiteness and Artinianness results. In the present section we study the vanishing
and Artinianness conditions of formal local cohomology modules as
our main result. Not so much is known about the mentioned properties. In \cite {A-D} Asgharzadeh and Divani-Aazar have investigated some properties of formal local cohomology modules. For instance they showed  that $ \mathfrak{F}^d_{\fa}(M)$ is Artinian for $d:= \dim M$. Here we give an alternative proof of it with more information on the attached primes of $ \mathfrak{F}^d_{\fa}(M)$:

\begin{proposition} \label{2.1} Let $\fa$ be an ideal of  a  local ring $(R,\fm)$ and $M$ a finitely generated $R$-module of dimension $d$. Then $ \mathfrak{F}^d_{\fa}(M)$ is Artinian. Furthermore
$$\begin{array}{ll} \Att_R (\mathfrak{F}^d_{\fa}(M)) =  \{\fp \in \Ass M: \dim R/\fp =d \} \cap V(\fa).
\end{array}
$$
\end{proposition}

{\bf Proof.}
By Independence Theorem we may assume that $\Ann M=0$ and so $d=\dim R$. As $H^d_{\fm} (M/\fa^n M)$ is right exact ($n \in \mathbb {N}$), we have
$$\begin{array}{ll} H^d_{\fm} (M/\fa^n M) & \cong H^d_{\fm} (R) \otimes_R M/\fa^n M
\\& \cong H^d_{\fm} (M) \otimes_R R/\fa^n
\\& \cong H^d_{\fm} (M)/\fa^n H^d_{\fm} (M).
\end{array}$$
Since $H^d_{\fm} (M)$ is an Artinian module so there exists an integer $n_0$ such that for all integer  $t \geq n_0$ we have  $\fa^t H^d_{\fm} (M) = \fa^{n_0} H^d_{\fm} (M)$. Then one can see that
$$\begin{array}{ll} \mathfrak{F}^d_{\fa}(M) \cong  H^d_{\fm} (M)/\fa^{n_0} H^d_{\fm} (M),
\end{array}$$
which is an Artinian module. By virtue of above equations and  \cite  [ Theorem 7.3.2]{Br-Sh}, the second claim is clear.
$\ \ \ \ \ \ \ \Box$

\begin{lemma} \label{2.2} Let $(R,\fm)$ be a complete local ring and $M$ a finitely generated $R$-module. Then $\Supp (\mathfrak{F}^0_{\fa}(M))=\bigcup_{\fp \in \Ass_R \mathfrak{F}^0_{\fa}(M)} V(\fp)$. Moreover $\Supp (\mathfrak{F}^0_{\fa}(M)) \cap V(\fa) \subseteq V(\fm)$ .
\end{lemma}

{\bf Proof.} To prove the claim, it is enough to consider that $\Ass_R \mathfrak{F}^0_{\fa}(M)=\{ \fp \in \Ass_R M: \dim R/(\fa +\fp)=0 \}$ (cf. \cite  [ lemma 4.1]{Sch}) .
$\ \ \ \ \ \ \ \Box$

Using  Lemma \ref{2.2} we are now able to prove Theorem \ref{1.1}:

 \textbf{Proof of Theorem \ref{1.1}.}

$(1)\Rightarrow (3)$ and  $(2)\Rightarrow (1)$ are obvious.

 $(3)\Rightarrow (2):$ By passing to the completion, we may assume that $R$ is complete (cf. \cite  [ Proposition 3.3]{Sch}).

We argue by induction on $t$. When $\ t=1$, there is nothing to prove, since Lemma \ref{2.2} and the assumptions imply that
$$\begin{array}{ll} \Supp (\mathfrak{F}^0_{\fa}(M)) = \Supp (\mathfrak{F}^0_{\fa}(M)) \cap V(\fa) \subseteq V(\fm).
\end{array}
$$
Hence $\ \mathfrak{F}^0_{\fa}(M)\ $ is Artinian. To this end note that $  \mathfrak{F}^0_{\fa}(M) $ is a finitely generated submodule of $M$. So suppose that $t> 1$ and the result has been proved for smaller values of $t$. Put $\overline {M} =M/H^0_{\fa}(M)$. From the short exact sequences
$$\begin{array}{ll} 0 \longrightarrow H^{0}_{\fa}(M ) \longrightarrow M \longrightarrow \overline{M} \longrightarrow 0
\end{array}
$$
and by \cite[Proposition 3.11]{Sch}), we get the following long exact sequence
$$\begin{array}{ll} ... \longrightarrow \mathfrak{F}^i_{\fa}(H^{0}_{\fa}(M )) \longrightarrow \mathfrak{F}^i_{\fa}(M) \longrightarrow \mathfrak{F}^i_{\fa}(\overline{M}) \longrightarrow \mathfrak{F}^{i+1}_{\fa}(H^{0}_{\fa}(M ))\longrightarrow ....
\end{array}
$$
As $\mathfrak{F}^j_{\fa}(H^{0}_{\fa}(M ))= H^j_{\fm}(H^0_{\fa}(M))$ is an Artinian $R$-module for every $j \in \mathbb {Z}$ (\cite  [ Theorem 7.1.3]{Br-Sh}) then, one can see that $\ \Supp (\mathfrak{F}^i_{\fa}(\overline {M})) \subseteq V(\fa)$ for all $i< t$. Hence, it is enough to show that $\mathfrak{F}^i_{\fa}(\overline {M})$ is Artinian, so we may assume that $H^0_{\fa}(M) = 0 $. Thus, there exists an $M$-regular element $x$ in ${\fa}$ such that from the short exact sequence
$$\begin{array}{ll} 0 \longrightarrow M \stackrel{x}{\rightarrow}  M \longrightarrow  M/xM= \widetilde{M}\longrightarrow 0
\end{array}
$$
we deduce the next long exact sequence
$$\begin{array}{ll} ... \longrightarrow \mathfrak{F}^i_{\fa}(M ) \stackrel{x}{\rightarrow} \mathfrak{F}^i_{\fa}(M) \longrightarrow \mathfrak{F}^i_{\fa}(\widetilde{M}) \longrightarrow \mathfrak{F}^{i+1}_{\fa}(M)\longrightarrow ....\ \ \ \ \ \ \ \ \ \ (\ast)
\end{array}
$$
Since   $\ \Supp (\mathfrak{F}^i_{\fa}(M)) \subseteq V(\fa) \ $ for all $i < t$, it follows from the above long exact sequence that $ \ \Supp (\mathfrak{F}^i_{\fa}(\widetilde{M})) \subseteq V(\fa) \ $ for all $\ i < {t-1} \ $. Hence, by induction hypothesis we have
$\ \mathfrak{F}^i_{\fa}(\widetilde{M}) \ $ is Artinian for all  $\ i < {t-1}$. Therefore in the view of $(\ast)$,    $\ ( 0 :_{\mathfrak{F}^i_{\fa}(M)}   x ) \ $ is  Artinian for all $\ i < t$ .

On the other hand since $\ \Supp (\mathfrak{F}^i_{\fa}(M)) \subseteq V(\fa)\ $ for all $i< t \ $, one can see that
$$\begin{array}{ll} \ \mathfrak{F}^i_{\fa}(M) = \bigcup ( 0 :_{\mathfrak{F}^i_{\fa}(M)}   {\fa^\alpha} ) \subseteq \bigcup ( 0 :_{\mathfrak{F}^i_{\fa}(M)} {x^\alpha}) \subseteq \mathfrak{F}^i_{\fa}(M)
\end{array}$$
so $ \ \mathfrak{F}^i_{\fa}(M) = \bigcup ( 0 :_{\mathfrak{F}^i_{\fa}(M)} {x^\alpha}) \ $. Therefore by \cite[Theorem 1.3]{Mel}, $\ \mathfrak{F}^i_{\fa}(M)\ $ will be Artinian for all $i<t$.

 $(2)\Rightarrow (4):$ Since $\mathfrak{F}^i_{\fa}(M)$ is $\fa$-adically complete for every $i \in \mathbb {Z}\ $(cf. \cite[Theorem 3.9]{Sch} or  \cite[Remark 3.1]{H-St} ), we get $\bigcap_{n}\fa ^n \mathfrak{F}^i_{\fa}(M)=0$. Moreover for all $i <t$,  $\ \mathfrak{F}^i_{\fa}(M)$ is Artinian. Hence,   there is an integer $u$ such that $\fa ^u \mathfrak{F}^i_{\fa}(M)=0$.

 $(4)\Rightarrow (3)$ is obvious.

$(1)\Rightarrow (5):$  By passing to the completion we may assume that $R$ is complete. We use induction on $t$. Let $t =1$. As $\Supp (\mathfrak{F}^0_{\fa}(M)) \subseteq V(\fm) $ so $\mathfrak{F}^0_{\fa}(M)  $ must be zero. Otherwise since
$$\begin{array}{ll} \ \emptyset \neq \Ass ( \mathfrak{F}^0_{\fa}(M) ) \subseteq \Supp ( \mathfrak{F}^0_{\fa}(M) ) \subseteq V({\fm})
\end{array}$$
then,
$$\begin{array}{ll} \ \fm \in \Ass ( \mathfrak{F}^0_{\fa}(M) )= \{ \fp \in \Ass M ;  \dim(R/\fa+\fp) = 0 \},
\end{array}$$
this is contradiction to $\depth M > 0$.

Now suppose that $ \depth M \geq t> 1$ and that  the result has been proved for smaller values of $t$. By this inductive assumption,   $\mathfrak{F}^i_{\fa}(M)=0$ for $i=0,1,...,t-2$ and it only remains for us to prove that $\mathfrak{F}^{t-1}_{\fa}(M)=0$.

Since $ \depth M > 1$  then, there exists $x \in {\fm} \ $ that is an $M$-regular element. Consider the short exact sequence
$$\begin{array}{ll} \ 0 \rightarrow M \stackrel{x^l}{\rightarrow} M \rightarrow M/x^l M = \bar{M} \rightarrow 0
\end{array}$$
for every $l$.  Thus, we have the following long exact sequence
$$\begin{array}{ll} \ ... \rightarrow \mathfrak{F}^{i-1}_{\fa}(\bar{M}) \rightarrow \mathfrak{F}^i_{\fa}(M) \stackrel{x^l}{\rightarrow} \mathfrak{F}^i_{\fa}(M)\rightarrow \mathfrak{F}^i_{\fa}(\bar{M}) \rightarrow ...
\end{array}$$
for every $l$.

As $\ \depth \bar{M}= \depth M -1> 0\ $  and for all $ i< t-1$,  $\ \Supp (\mathfrak{F}^i_{\fa}(\bar{M})) \subseteq V(\fm)$ then, by inductive hypothesis $\ \mathfrak{F}^i_{\fa}(\bar{M}) = 0 \ $  for all $\  i< t-1$ . Thus, for every $l$, $\ ( 0:_{\mathfrak{F}^{t-1}_{\fa}(M)} x^l)$ is a homomorphic image of $\mathfrak{F}^{t-2}_{\fa}(\bar{M})$. Hence, $\ ( 0:_{\mathfrak{F}^{t-1}_{\fa}(M)} x^l) = 0 \ $ for every $l$.

Take into account that by assumption $\Supp (\mathfrak{F}^i_{\fa}(M)) \subseteq V(\fm)$ for every $i<t$. Then,  $\ {\mathfrak{F}^{t-1}_{\fa}(M)}  = \cup ( 0:_{\mathfrak{F}^{t-1}_{\fa}(M)} x^l) = 0 \ $ . This completes the proof.
$\ \ \ \ \ \ \ \Box$

\section{Cosupport}

In this section we examine the cosupport of formal local cohomology.
The notion of cosupport was introduced by S. Yassemi in \cite{Y}. He defined the $\CoSupp_R M$ as the set of prime ideals $\fq$ such that there exists a cocyclic
homomorphic image $L$ of $M$ with $\fp \supseteq \Ann(L)$. His definition is equivalent to Melkersson-Schenzel's  definition for Artinian $R$-modules. Melkersson-Schenzel's definition of colocalization does not map Artinian $R$-module to Artinian $S^{-1}R$-module through colocalization at a multiplicative closed subset of $R$ (cf. \cite{Mel-Sch}). In this note we use the concept  of cosupport   has been introduced by A. Richardson \cite{R}. It maps Artinian $R$-modules to Artinian $S^{-1}R$-modules (when $R$ is complete). Also it is  suitable to investigate  formal local cohomology modules.

\begin{definition} \label{3.1} (cf. \cite {R}) Let $R$ be a ring and $M$ an $R$-module.
\begin{enumerate}
\item[(1)]  Let $S$ be a  multiplicative closed subset of $R$ and $D_{R}(-):= \Hom_{R}(-,E_{R})$, where $E_{R}$ is the injective hull of $\oplus R/\fm$, the sum running over all maximal ideals $\fm$ of $R$. The colocalization of $M$ relative to $S$ is the $S^{-1}R$-module $S_{-1}M=D_{S^{-1}R}(S^{-1}D_R(M))$. If $S=R \setminus \fp\ $ for some prime ideal $\fp\in \Spec(R)$, we write $^{\fp}M$ for $S_{-1}M$.
\item[(2)]  The cosupport of $M$ is defined as follows
$$\begin{array}{ll} \ \CoSupp_R M := \{ \fp\in \Spec(R):\ ^{\fp}M \neq 0 \}.
\end{array}$$
\end{enumerate}
\end{definition}
For brevity we often write $\CoSupp M$
for $\CoSupp_R M$
 when there is no ambiguity about the ring $R$.

 Below we recall some  properties of cosupport:

\begin{lemma}  \label{3.2} (cf. \cite [Theorem 2.7]{R}) Let $R$ be a ring and $M$ an $R$-module.
\begin{enumerate}
\item[(1)] $\CoSupp M =\Supp D_{R}(M)$.
\item[(2)] If $M$ is finitely generated, then $\CoSupp M = V(\Ann M)\cap \max(R)$.
\item[(3)] $\CoSupp M =\emptyset $ if and only if $M=0$.
\item[(4)]  $\CoSupp M \subseteq V(\Ann M)$.
\item[(5)] If $\ 0 \rightarrow M^{'} \rightarrow M \rightarrow M^{''} \rightarrow 0$ is exact, then $\CoSupp M = \CoSupp M^{'} \cup \CoSupp M^{''}$.
\item[(6)] If $M$ is representable, then $\CoSupp M = V(\Ann M)$.
\end{enumerate}
\end{lemma}

\begin{proposition} \label{3.3} Let $R$ be a ring and $M$ and $N$ be  $R$-modules. Then the following statements are true:
\begin{enumerate}
\item[(1)] $\CoSupp (M)$ is stable under specialization, i.e.
$$\begin{array}{ll} \  \fp \in \Cosupp (M) , \fp \subseteq \fq \Rightarrow \fq \in \Cosupp (M).
\end{array}$$
 \item[(2)] Let $M$ be a finite module, then  $\CoSupp (M \otimes_R N)\subseteq \Supp M \cap \CoSupp N$.
\end{enumerate}
\end{proposition}

{\bf Proof.}
\begin{enumerate}
\item[(1)] Let $\fp \in \Cosupp (M)$, then by definition $D_{R_{\fp}} (D_R(M)_{\fp}) $ is nonzero and so is $D_R(M)_{\fp}$. As  $0 \neq D_R(M)_{\fp}= (D_R(M)_{\fq})_{\fp R_{\fq}}$, then $D_R(M)_{\fq} \neq 0$. It implies that $^{\fq}M \neq 0$.

\item[(2)] Use \cite[2.5]{R} to prove.
%\item[(2)] Let $\fp \in \CoSupp (M \otimes_R N)$, then $0 \neq {^{\fp} (M \otimes_R N)}=M_{\fp} \otimes_{R_{\fp}} {^{\fp} N}$, by \cite[2.5]{11}. So $M_{\fp} \neq 0$ and ${^{\fp} N} \neq 0$. Hence $\fp \in \Supp M \cap \CoSupp N$.
\end{enumerate}
$\ \ \ \ \ \ \ \Box$

\begin{lemma} \label{3.4} Let $\fa$ be an ideal of a ring $R$. Let $N$ be an Artinian $R$-module with $\Att_R (N)\subseteq V(\fa)$. Then $\CoSupp N \subseteq V(\fa)$.
\end{lemma}

{\bf Proof.} Since $N$ is an Artinian module then, the following descending chain
$$\begin{array}{ll} \   \fa N \supseteq \fa^2 N \supseteq ... \supseteq \fa^n N  \supseteq ...
\end{array}$$
of submodules of $N$ is stable, i.e. there exists an integer $k$ that $\fa^k N= \fa^{k+1} N$. As  $\Att_R (N/\fa^k N)=\Att_R (N) \cap V(\fa)$ (cf. \cite[Proposition 5.2]{Mel-Sch}) and   $\CoSupp (N/\fa^k N) \subseteq V(\fa)$ by virtue of Proposition \ref{3.3}, hence, by passing to $N/\fa^k N$ we may assume that  $\fa^k N=0$.

Let $\fp \in \CoSupp N$ then, $^{\fp}N \neq 0$. Thus, for every $s \in S= R \setminus \fp$, $sN \neq 0$ (cf. \cite[2.1]{R}). On the other hand   $\bigcap_n \fa^nN=\fa^k N=0$, hence, for every $s \in S$, $sN \not \subseteq \fa^tN$. It follows that for all $s \in S$, $s \notin \fa^t$ and clearly  $\fp \in V(\fa)$.
$\ \ \ \ \ \ \ \Box$

\begin{corollary} \label{3.5} Let $i \in \mathbb {Z}$. Let $(R,\fm)$ be a local ring and  $M$ be a finitely generated $R$-module. Assume that  $\mathfrak{F}^i_{\fa}(M)$ is an Artinian $R$-module, then $\CoSupp \mathfrak{F}^i_{\fa}(M) \subseteq V(\fa)$.
\end{corollary}

{\bf Proof.} As $\mathfrak{F}^i_{\fa}(M)$ is Artinian and $\fa$-adically complete so, there exists an integer $k$ such that $\bigcap_{n \geq 1} \fa^n \mathfrak{F}^i_{\fa}(M)= \fa^k \mathfrak{F}^i_{\fa}(M)=0$. Hence, \cite[Proposition 7.2.11]{Br-Sh} implies that $\Att \mathfrak{F}^i_{\fa}(M)\subseteq V(\fa)$ and in the light of Lemma \ref{3.4} $\ \CoSupp \mathfrak{F}^i_{\fa}(M) \subseteq V(\fa)$.
$\ \ \ \ \ \ \ \Box$

\begin{remark} \label{3.6} Converse of Corollary \ref{3.5} is not true in general.  Let $R=k[\left| x \right|]$ denote the formal power series ring over a field $k$. Put $\fa=(x)R$.  Then
$$\begin{array}{ll} \  \CoSupp \mathfrak{F}^0_{\fa}(R) =\Supp D_R(D_R (H^{1}_{\fa}(R)))= \Supp H^{1}_{\fa}(R)\subseteq V(\fa)
\end{array}$$
but $\mathfrak{F}^0_{\fa}(R)$ is not Artinian .
\end{remark}

We now turn our attention to prove Theorem \ref{1.2}. For this reason we give a preliminary lemma:

\begin{lemma} \label{3.7} Let $(R, \fm)$ be a $d$-dimensional  local ring. Let $M$ be a finitely generated $R$-module. Then
$$\begin{array}{ll} \  \mathfrak{F}^c_{\fa}( M)\cong \mathfrak{F}^c_{\fa}( R)\otimes_{R}M,
\end{array}$$
where $c:= \dim R/\fa$.
\end{lemma}

{\bf Proof.} At first note that by definition of inverse limit, $\mathfrak{F}^j_{\fa}(-)$ preserves finite direct sum, for every $j \in \mathbb {Z}$. Furthermore  $\mathfrak{F}^c_{\fa}(-)$ is a right exact functor (cf. \cite[Theorem 4.5]{Sch}).
Hence, by Watts' Theorem ( \cite[Theorem 3.33]{Ro})  the claim is  proved. $\ \ \ \ \ \ \ \Box$

Lemma \ref{3.7} declares that   $\mathfrak{F}^c_{\fa}( R)=0$ if and only if $\mathfrak{F}^c_{\fa}( M)=0$ for all finitely generated $R$-module $M$.

In order to prove Theorem \ref{1.2} we  utilize the useful consequence of Gruson's Theorem (see, e.g., \cite[Corollary 4.3]{V}) allows us to reduce to the case $M=R$ when considering the cosupport of top formal local cohomology modules:

\textbf{Proof of Theorem \ref{1.2}.}

 (1): Since $\mathfrak{F}^c_{\fa}( M) \cong \mathfrak{F}^c_{\fa (R/J)}( M)$,  by Independence Theorem \cite [4.2.1]{Br-Sh}, we may replace $R$ by $R/J$ to assume that $M$ is faithful. Note that for $\dim R/(\fa, J) < c$, there is nothing to prove because, $\mathfrak{F}^c_{\fa}( M)= 0$.

In the view of lemma \ref{3.7} and \cite[Proposition 2.5]{R}, for every $\fp \in \Spec R$
$$\begin{array}{ll} \    ^{\fp}\mathfrak{F}^c_{\fa}( M) \cong M_{\fp} \otimes_{R_{\fp}}\  ^{\fp}\mathfrak{F}^c_{\fa}( R).
\end{array}$$

As  $M_{\fp}$ is a faithful $R_{\fp}$-module, \cite[Corollary 4.3]{V} implies that $M_{\fp} \otimes_{R}\  ^{\fp}\mathfrak{F}^c_{\fa}( R)=0$ if and only if
 $ ^{\fp}\mathfrak{F}^c_{\fa}( R)=0$, which completes the proof.

(2): To prove,  we use the localization instead of colocalization in the proof of $(1)$.
$\ \ \ \ \ \ \ \Box$

\section{Coassociated primes}

Let $M$ be an $R$-module. A prime ideal $\fp$ of $R$ is called
a coassociated prime of $M$ if there exists a cocyclic homomorphic image
$L$ of $M$ such that $\fp = \Ann(L)$. The set of coassociated prime ideals of
$M$ is denoted by $\Coass_R(M)$ (cf. \cite{Y}). When the ambient
R is understood, we will often write $\Coass(M)$ instead of $\Coass_R(M)$.

Note that for an Artinian module the set of coassociated primes is finite. In this section $(R,\fm)$ is a local ring and we denote by $D_R(M)=\Hom_R (M, E(R/\fm)) $ the Matlis dual of $R$-module $M$, where $E(R/\fm)$ is the injective hull of residue field,  so in this case $\Coass(M)=\Ass D_R(M)$.

Among other results, we will see that under certain assumptions $ \Cosupp_R (\mathfrak{F}^i_{\fa}( M))$  as a subset of $\Spec R$ is  closed in the Zariski topology for some $i \in \mathbb {Z}$.

\begin{lemma} \label{4.1} Let $(R, \fm)$ be a  local ring and $M$ an $R$-module. Then the following statements are true:
\begin{enumerate}
\item[(1)] $\Coass (M) \subseteq \CoSupp (M)$.
\item[(2)] Every minimal element of $\CoSupp (M)$ belongs to $\Coass (M)$.
\item[(3)] For any Noetherian $\widehat{R}$-module $M$,   $\Coass (M) = \CoSupp (M) \subseteq \{ \fm \}$, where $\widehat{R}$ denotes the  $\fm$-adic completion of $R$.
\end{enumerate}
\end{lemma}

{\bf Proof.}
\begin{enumerate}
\item[(1)] Let $\fp \in \Coass (M)$, then, it implies that $0\neq \Hom_{R_{\fp}}(R_{\fp}/\fp R_{\fp},D_{R}(M)_{\fp})$. Note that  it remains nonzero by taking $\Hom_{R_{\fp}}(-,E_{R_{\fp}}(R_{\fp}/\fp R_{\fp}))$ and consequently $\fp \in \CoSupp (M)$.

\item[(2)]   Let $\fp \in \min \CoSupp (M)=  \min \Supp D_{R}(M)$, so, $\fp \in \min \Ass D_{R}(M)$. It follows that $\fp \in \min \Coass (M)$.

\item[(3)] It is clear by (1) and (2).$\ \ \ \ \ \ \ \Box$
\end{enumerate}

It should be noted that $\Supp (\mathfrak{F}^i_{\fa}( M))$ is closed when $\Ass (\mathfrak{F}^i_{\fa}( M))$ is finite. In fact for a local Gorenstein ring $(R,\fm)$, $\ \Ass (\mathfrak{F}^i_{\fa}( R))= \Ass D_{R}(H^{\dim R-i}_{\fa}(R))$ see \cite{H} for details. Take into account that it is not finite in general ( see \cite{H} or \cite [Remark 2.8(vi)]{A-D}).

\begin{lemma} \label{4.2} Let $(R, \fm)$ be a local ring and $M$ be an $R$-module. The set of  minimal primes in $\CoSupp(M)$ is finite if and only if $\CoSupp(M)$ is a closed subset of $\Spec R$.
\end{lemma}

{\bf Proof.} Let $\CoSupp (M)= V(\fb)$ for some ideal $\fb$ of $R$. As $R$ is Noetherian then so is $R/\fb$. It turns out that the set of minimal elements of $\CoSupp(M)$ is finite.

For the reverse direction, let $\fp_{1},..., \fp_{t}$ be the minimal prime ideals of $\CoSupp (M)$. Put $\fq:= \cap_{i} \fp_{i}$. We claim that $\CoSupp  M= V(\fq)$.

It is clear that $\CoSupp (M)\subseteq V(\fq)$. For the opposite direction assume that there is a prime ideal $Q \supset \fq$. Then $Q \supset \fp_{j}$, for some $1 \leq j \leq t$ so the proof follows by \ref{3.3}(1). $\ \ \ \ \ \ \ \Box$

We deduce from above  lemma that the cosupport of formal local cohomology module is closed, whenever its set of coassociated primes  is finite.  Therefore if one of the situations in Theorem \ref{1.1} is true,  the cosupport of formal local cohomology module is closed. Also  $\CoSupp  \mathfrak{F}^{\dim M}_{\fa}(M)$ is closed, as $\mathfrak{F}^{\dim M}_{\fa}(M)$ is Artinian, whenever $M$ is a finitely generated module over a local ring $(R, \fm)$ (cf. \cite[Lemma 2.2]{A-D}).

Take into account that when $R$ is a complete local Gorenstein ring and $\mathfrak{F}^i_{\fa}(M)$ is assumed to be either Noetherian or Artinian module, then
$$\begin{array}{ll} \   \Cosupp ( \mathfrak{F}^i_{\fa}(M))= \Supp   H^{\dim R-i}_{\fa}(M,R).
\end{array}$$

By virtue of \cite[Theorem 2.7]{A-D}, for a Cohen-Macaulay ring $R$ with $\Ht \fa > 0$,
$\mathfrak{F}^{\dim R/\fa}_{\fa}(R)$
 is not Artinian. Moreover $\mathfrak{F}^{\dim M/\fa M}_{\fa}(M)$
 is not finitely generated for
$\dim M/\fa M > 0$ (cf. \cite[Theorem 2.6(ii)]{A-D}). Below we give an alternative proof:

\begin{theorem} \label{4.3} Let $\fa$ be an ideal of a local ring $(R,\fm)$ and $M$ a finitely generated R-module. Assume that $\dim M/\fa M > 0$. Then $\mathfrak{F}^{\dim M/\fa M}_{\fa}( M)$ is not a finitely generated $R$-module.
\end{theorem}

{\bf Proof.} Put $c:= \dim M/\fa M$. In the contrary assume that $\mathfrak{F}^{c}_{\fa}( M)$ is a finitely generated $R$-module. Let $x \in \fm$ be a parameter element of $M/\fa M$. Hence, \cite[Theorem 3.15]{Sch} implies the following long exact sequence
$$\begin{array}{ll} \  ... \rightarrow \Hom (R_x, \mathfrak{F}^{c}_{\fa}( M)) \rightarrow \mathfrak{F}^{c}_{\fa}( M) \rightarrow \mathfrak{F}^{c}_{(\fa,x)}( M) \rightarrow ...,
\end{array}$$
where $i \in \mathbb {Z}$.
As $\dim M/(\fa,x) M < \dim M/\fa M\ $ then, $\mathfrak{F}^{c}_{(\fa,x)}( M)=0$. Now let $f \in \Hom (R_x, \mathfrak{F}^{c}_{\fa}( M))$. Fix an arbitrary integer $n$, so
$$\begin{array}{ll} \  f(1/x^n)=x^m f(1/x^{m+n}) \in x^m \mathfrak{F}^{c}_{\fa}( M),
\end{array}$$
for every integer $m$. It implies that $f(1/x^n) \in \bigcap_{m} x^m \mathfrak{F}^{c}_{\fa}( M)=0$ by Krull's Theorem and hence, $f=0$. Now it follows that $\mathfrak{F}^{c}_{\fa}( M)=0$, which is a contradiction, see \cite[Theorem 4.5]{Sch}.
$\ \ \ \ \ \ \ \Box$

Now we examine the set of coassociated primes of top formal local cohomology to show that by some  assumptions on $R$, it could be finite.

\begin{proposition} \label{4.4} Let $\fa$ be an ideal of a  complete Gorenstein local  ring $(R, \fm)$ and $c:= \dim R/\fa$. Let $M$  be a  finitely generated  $R$-module. Then
$$\begin{array}{ll} \   \Coass_R (\mathfrak{F}^c_{\fa}( M)) = \Supp_R (M)  \cap \Ass (H^{\Ht \fa}_{\fa}(R)).
\end{array}$$
In particular,  $\ \Coass_R (\mathfrak{F}^c_{\fa}( M))$ is finite.
\end{proposition}

{\bf Proof.}
$$\begin{array}{ll} \   \Coass_R (\mathfrak{F}^c_{\fa}( M))& = \Coass_R (\mathfrak{F}^c_{\fa}( R)\otimes_R M)
\\& = \Supp_R M \cap  \Coass_R (\mathfrak{F}^c_{\fa}( R))
\\& = \Supp_R M \cap \Ass_R (H^{\Ht \fa}_{\fa}(R))
\end{array}$$
where the first equality is clear by Lemma \ref{3.7}, the second equality follows by \cite[Theorem 1.21]{Y}.
$\ \ \ \ \ \ \ \Box$

It should be noted that by hypotheses in Proposition \ref{4.4}, $\Ht \fa = \grade_R \fa$  and it is well-known
that $\Ass_R (H^{\grade_R \fa}_{\fa}(R))$
 is finite, cf. \cite{Br-L}.

\begin{corollary} \label{4.5} Keep the notations and hypotheses in Proposition \ref{4.4},
$$\begin{array}{ll} \     \mathfrak{F}^c_{\fa}( M)= 0 \Longleftrightarrow  \Supp_R (M)  \cap \Ass (H^{\Ht \fa}_{\fa}(R)) =\emptyset.
\end{array}$$
\end{corollary}

%{\bf Proof.} It is an immediate consequence of  (4.3).$\ \ \ \ \ \ \ \Box$
%{\bf Proof.} Consider $\mathfrak{F}^d_{\fa}( M)= 0 \ $ then $\Coass (\mathfrak{F}^d_{\fa}( M))= \emptyset$, then apply (4.3).

%For the converse, by 3.2(3) it is enough to show that $\CoSupp (\mathfrak{F}^d_{\fa}( M)) = \emptyset$ and it is clear by virtue of  (4.1)(2).
%$\ \ \ \ \ \ \ \Box$

\begin{proposition} \label{4.6} Let $i \in \mathbb {Z}$. Let $\fa \subset R$ be an ideal of a ring $R$. If $\Coass_R \mathfrak{F}^i_{\fa}( R)$ is finite, then so is  $\Coass_R \mathfrak{F}^i_{\fa}( R/H^0_{\fa}(R))$.
\end{proposition}

{\bf Proof.} Consider the exact sequence
$$\begin{array}{ll} \    0 \rightarrow H^0_{\fa}(R) \rightarrow R \rightarrow  R/H^0_{\fa}(R)=\overline{R} \rightarrow 0.
\end{array}$$
It provides the following long exact sequence
$$\begin{array}{ll} \    ... \rightarrow \mathfrak{F}^i_{\fa}(H^0_{\fa}(R)) {\stackrel{\psi}{\rightarrow}} \mathfrak{F}^i_{\fa}(R) \stackrel{\varphi}{\rightarrow} \mathfrak{F}^i_{\fa}(\overline{R}) \rightarrow \mathfrak{F}^{i+1}_{\fa}(H^0_{\fa}(R)) \rightarrow ..., \ \ \ \ \ \ \ \ (\ast)
\end{array}$$
for every $i$.

As $\mathfrak{F}^i_{\fa}(H^0_{\fa}(R))=H^i_{\fm}(H^0_{\fa}(R))$ is Artinian, it follows that $\Coass (\mathfrak{F}^i_{\fa}(H^0_{\fa}(R)))$ is finite.

By virtue of $(\ast)$, we get the following short exact sequence
$$\begin{array}{ll} \    0 \rightarrow  U \rightarrow  \mathfrak{F}^i_{\fa}(\overline{R}) \rightarrow U^{'} \rightarrow 0,
\end{array}$$
where $U=\coker \psi$ and $U^{'}=\coker \varphi$. It implies that $\Coass \mathfrak{F}^i_{\fa}( \overline{R})$ is finite. To this end note that $\Coass (U)$ is finite by \cite[Theorem 1.10]{Y} and $\Coass (U^{'})$ is finite as $\mathfrak{F}^{i+1}_{\fa}(H^0_{\fa}(R))$
is Artinian.
$\ \ \ \ \ \ \ \Box$

Now we are going to give more information on the last non-vanishing formal
local cohomology module.

\begin{theorem} \label{4.7} Let $(R,\fm)$ be a local ring of dimension $d >1$. Let $\mathfrak{F}^{d}_{\fa}(R)=0$. Then:
\begin{enumerate}
\item[(1)] If $\fp \in \Coass \mathfrak{F}^{d-1}_{\fa}(R)$, then it implies that $\dim (R/(\fa,\fp))= d-1$.
\item[(2)] $\Assh(R) \cap \Coass \mathfrak{F}^{d-1}_{\fa}(R)\subseteq \{ \fp \in \Spec R: \dim R/\fp=d, \Rad (\fa+ \fp) \neq \fm \}$.
\item[(3)] If  $ \Coass \mathfrak{F}^{d-1}_{\fa}(R)\subseteq \Assh(R)$, then $\{ \fp \in \Spec R: \dim (R/(\fa,\fp))= d-1 \} \subseteq \Coass \mathfrak{F}^{d-1}_{\fa}(R)$.
\end{enumerate}
\end{theorem}

{\bf Proof.}
\begin{enumerate}
\item[(1)] Let $\fp \in \Coass \mathfrak{F}^{d-1}_{\fa}(R)$.  As $ \mathfrak{F}^{d}_{\fa}(R)=0$ then, by \cite[Theorem 4.5]{Sch} we have
$$\begin{array}{ll} \  \dim R/(\fa,\fp) \leq \dim R/\fa \leq d-1.
\end{array}$$
On the other hand $\fp \in \Coass (R/\fp \otimes_R \mathfrak{F}^{d-1}_{\fa}(R))$, because
$$\begin{array}{ll} \  \Coass (R/\fp \otimes_R \mathfrak{F}^{d-1}_{\fa}(R))= \Supp R/\fp \cap \Coass \mathfrak{F}^{d-1}_{\fa}(R).
\end{array}$$
 It yields with the similar argument to lemma \ref{3.7} that $0 \neq R/\fp \otimes_R \mathfrak{F}^{d-1}_{\fa}(R)=\mathfrak{F}^{d-1}_{\fa}(R/\fp)$. So, we have $ \dim R/(\fa,\fp)\geq d-1$. It completes the proof.

\item[(2)] Let $ \fp \in \Assh(R) \cap \Coass \mathfrak{F}^{d-1}_{\fa}(R)$. Then, similar to $(1)$,  $\mathfrak{F}^{d-1}_{\fa}(R/\fp)\neq 0$ and  moreover $\Rad (\fa+ \fp) \neq \fm$. To this end note that if $\Rad (\fa+ \fp) = \fm$ then,  $ \mathfrak{F}^{d-1}_{\fa}(R/\fp)=0 $ by Grothendieck's vanishing Theorem.

\item[(3)] Let $\fp \in \Spec R$ and  $\dim (R/(\fa,\fp))= d-1 $. Then, it follows that $\emptyset \neq \Coass \mathfrak{F}^{d-1}_{\fa}(R/\fp) = \Supp(R/\fp) \cap \Coass\mathfrak{F}^{d-1}_{\fa}(R)$.
Let $\fq \in \Coass\mathfrak{F}^{d-1}_{\fa}(R)$ then, $\fq \supseteq \fp$, but by assumption $\fq$ is minimal so we deduce that $\fq=\fp$.
$\ \ \ \ \ \ \ \Box$
\end{enumerate}

\begin{remark} \label{4.8}
 The inclusion in Theorem \ref{4.7}(2) is not an equality in general. For example Let $R=k[[x,y,z]]$ denote the formal power series ring in three  variables over a field $k$. Let  $\fa=(x,y)$ be an ideal of $R$ which is of  dimension one and put $\fp=0$. It is clear that  $\mathfrak{F}^{3-1}_{\fa}(R)=0=\mathfrak{F}^{3}_{\fa}(R)$, that is $\Coass \mathfrak{F}^{3-1}_{\fa}(R)= \emptyset$.
\end{remark}

\begin{lemma} \label{4.9} Let $(R,\fm)$ be a complete local ring and $\fa$ an ideal of $R$. Let $\fp$ be a minimal prime ideal of $\fa$. Then   $\fq \in \Coass_R(\widehat{R_{\fp}})$ implies that $\fq \subseteq \fp$, where the functor $\ \widehat{.}\ $ denotes the completion functor.
\end{lemma}

{\bf Proof.} The proof is straightforward. Let $\fq \in \Coass_R(\widehat{R_{\fp}})$, then
$$\begin{array}{ll} \   0 \neq  \Hom_R (R/\fq , \Hom_R( \widehat{R_{\fp}}, E_R(R/\fm))) = \Hom_R (R/\fq \otimes_R \widehat{R_{\fp}}, E_R(R/\fm)).
\end{array}$$
It yields that
$$\begin{array}{ll} \   0 \neq  R/\fq \otimes_R \widehat{R_{\fp}} = R/\fq \otimes_R R_{\fp} \otimes_{R_{\fp}} \widehat{R_{\fp}}.
\end{array}$$
It is clear that $R_{\fp}/\fq R_{\fp} \neq 0$ and so $\fq$ must be contained in $\fp$.
$\ \ \ \ \ \ \ \Box$

\textbf{Proof of Theorem \ref{1.4}}

For $d>2$ and $d=1$, the claim is clear.

 Let $d=2$. Suppose that $\fp_1,..., \fp_r$ are the minimal prime ideals of $\fa$. Put $S= \bigcap^r_{i=1}(R \setminus \fp_i)$ and choose $y \in \fm \setminus \bigcup^r_{i=1} \fp_i$. By  \cite[Theorem 2.2.4]{Br-Sh}, for any $n \in \mathbb {N}$ we have
 $$\begin{array}{ll} \ 0 \rightarrow H^0_{\fm}(R/\fa^n) \rightarrow R/\fa^n \rightarrow D_{(y)}(R/\fa^n) \rightarrow H^1_{\fm}(R/\fa^n) \rightarrow 0,
\end{array}$$
where $D_{(y)}(R/\fa^n)$ is the $(y)$-transform functor. One can see that $D_{(y)}(R/\fa^n)\cong R_S/\fa^n R_S$, so we get the following exact sequence
 $$\begin{array}{ll} \ 0 \rightarrow H^0_{\fm}(R/\fa^n) \rightarrow R/\fa^n \rightarrow R_S/\fa^n R_S \rightarrow H^1_{\fm}(R/\fa^n) \rightarrow 0.
\end{array}$$
Furthermore $R_S/\fa^n R_S \cong \oplus^r_{i=1} R_{\fp_i}/ \fa^n R_{\fp_i}$.
All the modules in the above exact sequence satisfy the Mittag-Leffler condition so by applying inverse limits we get
$$\begin{array}{ll} \ 0  \rightarrow R/\mathfrak{F}^{0}_{\fa}(R) \rightarrow \oplus^r_{i=1} \widehat{R_{\fp_i}} \rightarrow \mathfrak{F}^{1}_{\fa}(R) \rightarrow 0.
\end{array}$$
It yields that $\Coass_R (\mathfrak{F}^{1}_{\fa}(R)) \subseteq \bigcup^r_{i=1} \Coass_R(\widehat{R_{\fp_i}}) \subseteq \Coass_R (\mathfrak{F}^{1}_{\fa}(R)) \cup \{ \fm\}$. In the view of lemma \ref{4.9}, $\Coass_R (\mathfrak{F}^{1}_{\fa}(R)) = \bigcup^r_{i=1} \Coass_R(\widehat{R_{\fp_i}})$. Now the claim is proved by \cite[Beispiel 2.4]{Z}. To this end note that $\Coass_R(\widehat{R_{\fp_i}})=\Coass_{R_{\fp_i}}(\widehat{R_{\fp_i}}) \cap R$ for every $i \in \{ 1,...,r \}$.
$\ \ \ \ \ \ \ \Box$

\begin{remark} \label{4.10} Keep the notations and hypotheses in Theorem \ref{1.4}  and let $M$ be a finitely generated $R$-module. As $R$ is complete so by Cohen's structure Theorem, there exists a Gorenstein local ring $(S,\fn)$ where $R$ is a homomorphic image of  $S$ and $\dim R= \dim S$. Then by virtue of \ref{3.7} we have
 $$\begin{array}{ll} \ \Ass_R H^1_{\fa S}(M,S) \subseteq \Coass \mathfrak{F}^{d-1}_{\fa}(R)
\end{array}$$
is finite.
\end{remark}

\begin{acknowledgement} My thanks are due to my phd. adviser, Professor Peter Schenzel, for his guidance to prepare this paper and useful hints and to the reviewer for suggesting several improvements. Some parts of this paper  was written while the author was at Oberwolfach: Representations of Finite Groups, Local Cohomology and Support. Many thanks to the organisers.
\end{acknowledgement}

%%%%%%%%%%%%%%%%%%%%%%%%%%%%%%%%%%%%%%%%%%%%%%%%%%%


\begin{thebibliography}{99}

\bibitem[1]{A-D}{M. Asgharzadeh  and K. Divaani-Aazar}, {\it finiteness properties of formal local cohomology modules and Cohen-Macaulayness},    Commun. Alg. {\bf 39}, no. 3, (2011), 1082-1103(22).

\bibitem[2]{Br-Sh}{M. Brodmann and R.Y. Sharp}, {\it Local cohomology:
an algebraic introduction with geometric applications}, Cambridge Univ.
Press, {\bf 60}, Cambridge, (1998).

\bibitem[3]{Br-L}{M. Brodmann and A. Lashgari  Faghani}, {\it A finiteness result for associated primes of local cohomology modules}, Proc. Amer. Math. Soc. {\bf 128}, no. 10, 2851-2853; (2000).
Press, {\bf 60}, Cambridge, (1998).

\bibitem[4]{H}{M. Hellus}, {\it Local Cohomology and Matlis Duality}, Habilitationsschrift, Leipzig, (2006).



\bibitem[5]{H-St}{M. Hellus and J. St\"{u}ckrad}, {\it On endomorphism rings of local cohomology modules}, Proc. Amer. Math. Soc. {\bf 136} , no. 7, 2333-2341; (2008).



\bibitem[6]{Hu}{C. Huneke}, {\it Problems on local cohomology: Free resolutions in commutative algebra and algebraic geometry}, (Sundance, UT, 1990), 93-108, Jones and Bartlett, (1992).

\bibitem[7]{Hu-K-Mar}{C. Huneke, D. Katz and T. Marley }, {\it  On the support of local cohomology}, J.  Algebra {\bf 322} (2009) 3194-3211.

\bibitem[8]{Ma}{A. Mafi}, {\it Results of formal local cohomology modules},  To appear in Bull. Malays. Math. Sci. Soc.

\bibitem[9]{Mel}{L. Melkersson}, {\it on asymptotic stability for sets of prime ideals connected with the powers of an ideal},  Math. Proc. Camb. Phil. Soc. {\bf 107}, (1990),  267-271.

\bibitem[10]{Mel-Sch}{L. Melkersson and P. Schenzel}, {\it  The co-localization of an Artinian module}, Proc. Edinburgh Math. Soc. {\bf 38}, 121-131 (1995).

\bibitem[11]{O}{A. Ogus}, {\it Local cohomological dimension of Algebraic Varieties}, Ann. of Math. , {\bf 98}(2), (1973), 327-365.


\bibitem[12]{P-S}{C. Peskine and L. Szpiro}, {\it Dimension projective finie et cohomologie locale}, Publ. Math. I.H.E.S., {\bf 42}, (1973), 323-395.

\bibitem[13]{R}{A. S. Richardson}, {\it Co-localization, co-support and local cohomology}, Rocky Mountain J. of Math., {\bf 36}, 5, (2006), 1679-1703.

\bibitem[14]{Ro}{J. Rotman}, {\it An Introduction to Homological Algebra}, Academic Press, Orlando, FL, 1979.


\bibitem[15]{Sch}{P. Schenzel},  {\it On formal local cohomology and
connectedness}, J. Algebra, {\bf 315}(2), (2007), 894-923.


\bibitem[16]{V}{W. Vasconcelos},  {\it Divisor Theory in Module Categories }, North-Holland Math. Stud., vol. {\bf 14} , Notas de Matemática
(Notes on Mathematics), vol. 53, North-Holland Publishing Co./American Elsevier Publishing Co., Inc., Amsterdam/Oxford/New York, (1974).

\bibitem[17]{Y}{S. Yassemi},  {\it Coassociated primes }, Comm. Algebra, {\bf 23}(4), 1473-1498 (1995).


\bibitem[18]{Z}{H. Z\"{o}schinger},  {\it Der Krullsche Durchschnittssatz f\"{u}r kleine Untermoduln}, Arch. Math. (Basel), {\bf 62}(4),
(1994), 292-299.




\end{thebibliography}
\end{document}